\documentclass[12pt]{amsart}
\usepackage{amsmath}
\usepackage{amssymb}
\usepackage{amsmath,amssymb,amsthm,amscd}
\usepackage{relsize}

\usepackage{txfonts}
\usepackage{amsbsy}
\usepackage{graphicx,color}

\numberwithin{equation}{section}

\newtheorem{prop}{Proposition}[section]
\newtheorem{theo}{Theorem}[section]
\newtheorem{lemm}{Lemma}[section]

\def\begeq{\begin{equation}}
\def\endeq{\end{equation}}

\begin{document}

\title{Sufficient and Necessary Condition of the $L_p$-Brunn-Minkowski Inequality Conjecture for $p\in[0,1)$}
\author{Shi-Zhong Du}
\thanks{The author is partially supported by NSFC (12171299), and GDNSF (2019A1515010605)}
  \address{The Department of Mathematics,
            Shantou University, Shantou, 515063, P. R. China.} \email{szdu@stu.edu.cn}

\renewcommand{\subjclassname}{%
  \textup{2010} Mathematics Subject Classification}
\date{Jan. 2022}

\begin{abstract}
   The $L_p-$Brunn-Minkowski inequality plays a central role in the Brunn-Minkowski theory proposed by Firey \cite{F} in 60's and developed by Lutwak \cite{L2,L3} in 90's, which generalizes the classical Brunn-Minkowski inequality by $L_p-$sum of convex bodies. The inequality has been established by Firey for $p>1$ and later been conjectured by B\"{o}r\"{o}czky-Lutwak-Yang-Zhang \cite{BLYZ} for $p\in[0,1)$. (see also \cite{KM,CHLL}) The validity of this conjecture was verified for the planar case in \cite{BLYZ}, and for the higher dimensional case when $p$ closing to $1$ by Chen-Huang-Li-Liu \cite{CHLL}(see also a local version by Kolesnikov-Milman \cite{KM}). In this short note, we give a simple argument clarifying the equivalence between the full conjecture and a lower bound of the third eigenvalue of the Aleksandrov's problem.
\end{abstract}

\maketitle\markboth{$L_p$-Brunn-Minkowski inequality}{$L_p$-Minkowski problem}

\tableofcontents

\section{Introduction}

In the classical Brunn-Minkowski theory, the Minkowski sum of two convex bodies $K, L$ is defined by
   $$
    (1-\lambda)K+\lambda L\equiv\Big\{(1-\lambda)x+\lambda y\Big|\ x\in K, y\in L\Big\}.
   $$
A famous fact
   \begin{equation}\label{e1.1}
     V^{\frac{1}{n}}((1-\lambda)K+\lambda L)\geq(1-\lambda)V^{\frac{1}{n}}(K)+\lambda V^{\frac{1}{n}}V(L), \ \ \forall\lambda\in[0,1]
   \end{equation}
was known to be the classical Brunn-Minkowski inequality with equality if and only if $K$ and $L$ are homothetic. This inequality is also equivalent to its log-concavity form
   \begin{equation}\label{e1.2}
     V((1-\lambda)K+\lambda L)\geq V^{1-\lambda}(K)V^\lambda(L), \ \ \forall\lambda\in[0,1].
   \end{equation}
In the 60's, Firey \cite{F} generalized Minkowski sum to its $L_p$-analogue by the so called $L_p-$Minkowski sum
  $$
   (1-\lambda)\cdot K+_p\lambda\cdot L\equiv S(h_{K,L,p,\lambda})
  $$
for $p>1$ and support functions $h_K, h_L$ of $K,L$ respectively, where $S(h_{K,L,p,\lambda})$ is the convex body determined by support function
   $$
    h_{K,L,p,\lambda}\equiv\Big((1-\lambda)h_K^p(x)+\lambda h^p_L(x)\Big)^{\frac{1}{p}}.
   $$
 He also established the Brunn-Minkowski-Firey inequality
   \begin{equation}\label{e1.3}
     V\Big((1-\lambda)\cdot K+_p\lambda\cdot L\Big)\geq V^{1-\lambda}(K)V^\lambda(L), \ \ \forall \lambda\in[0,1]
   \end{equation}
for $p>1$, with equality if and only if $K=L$. It is natural to look into the case $p\in[0,1)$ for $h_{K,L,p,\lambda}$ defined above in case $p\in(0,1)$ and for
   $$
    h_{K,L,0,\lambda}\equiv h_K^{1-\lambda}(x)h^\lambda_L(x)
   $$
in case $p=0$. However, when $p\in[0,1)$, the function $h_{K,L,p,\lambda}$ may not be any support function of a convex body. As suggested by B\"{o}r\"{o}czky-Lutwak-Yang-Zhang \cite{BLYZ}, one can still define the $L_p-$Minkowski sum by
   $$
    (1-\lambda)\cdot K+_p\lambda\cdot L\equiv W(h_{K,L,p,\lambda}),
   $$
where
   $$
     W(h_{K,L,p,\lambda})\equiv\bigcap_{x\in{\mathbb{S}}^n}\Big\{y\in{\mathbb{R}}^{n+1}\Big|\ x\cdot y\leq h_{K,L,p,\lambda}\Big\}
   $$
is the Wulff body determined by the largest convex body whose support function is no greater than the Wulff function $h_{K,L,p,\lambda}$. With the help of this extension, they proposed the following conjecture for $p\in[0,1)$. (see also Kolesnikov-Milman \cite{KM} and  Chen-Huang-Li-Liu \cite{CHLL})\\

\noindent\textbf{The conjecture of global $L_p-$Brunn-Minkowski inequality:} \eqref{e1.3} holds for $p\in[0,1)$ and all origin symmetric convex bodies $K, L$ in ${\mathbb{R}}^{n+1}$ with equality holds if and only if $K=L$.\\

In the case $p=0$, the inequality was also named log-Brunn-Minkowski inequality in the literatures. There are some simple examples constructed in \cite{BLYZ} showing that \eqref{e1.3} does not hold for non-origin-symmetric convex bodies. Moreover, the validity of the conjecture was verified by B\"{o}r\"{o}czky-Lutwak-Yang-Zhang \cite{BLYZ} for the planar case, and was verified by Chen-Huang-Li-Liu \cite{CHLL} for higher dimensional case when $p$ closing to one. (see also Kolesnikov-Milman \cite{KM} for a localized version of the inequality in the same range of $p$ closing to one) The main purpose of this note is to prove the following result.

\begin{theo}\label{t1.1}
  The $L_p$-Brunn-Minkowski inequality holds for original symmetric convex bodies $K$ and $p\in[0,1)$ if and only if the third eigenvalue $\lambda_3(K)$ of the Aleksandrov's eigenvalue problem
   \begin{equation}\label{a.1}
     -U^{ij}(\nabla^2_{ij}\phi+\phi\delta_{ij})=\lambda h_K^{-1}\det(\nabla^2h_K+h_KI)\phi
   \end{equation}
  for support function $h_K$ and $[U^{ij}]\equiv co[\nabla^2_{ij}h_K+h_K\delta_{ij}]$ satisfies that
    \begin{equation}\label{a.2}
      \lambda_3(K)\geq1.
    \end{equation}
\end{theo}

As pointed out by Professor Shi-Bing Chen, it is well known that the argument in the paper \cite{CHLL} actually shows that the sufficiency of \eqref{a.2} for the validity of the conjecture, although it is not claimed explicitly. So, in this short note, we just clarify their results together with the necessity.

There are several equivalent forms of global $L_p-$Brunn-Minkowski inequality. The first one was well known as the global $p-BM$-inequality
    \begin{equation}\label{e1.4}
      V((1-\lambda)\cdot K+_p\lambda\cdot L)\geq\Big((1-\lambda)V^{\frac{p}{n}}(K)+\lambda V^{\frac{p}{n}}(L)\Big)^{\frac{n}{p}}, \ \ \forall \lambda\in[0,1].
    \end{equation}
The second one was shown by B\"{o}r\"{o}czky-Lutwak-Yang-Zhang using the so called $L_p-$Minkowski inequality.

\begin{theo}\label{t1.2} (\cite{BLYZ}, page 1982, Lemma 3.1)
 For each $p>0$, when restricted to origin-symmetric convex bodies in ${{\mathbb{R}}^{n+1}}$, the $L_p-$Brunn-Minkowski inequality \eqref{e1.3} is equivalent to $L_p-$Minkowski inequality
    \begin{equation}\label{e1.5}
      \Bigg(\int_{{\mathbb{S}}^n}\Bigg(\frac{h_L}{h_K}\Bigg)^pd\overline{V}_K\Bigg)^{\frac{1}{p}}\geq\Bigg(\frac{V(L)}{V(K)}\Bigg)^{\frac{1}{n}}
    \end{equation}
 where $\overline{V}_K\equiv\frac{V_K}{V(K)}$ is the normalized cone-volume measure of $K$ and equality holds if and only if $K, L$ are dilates.
\end{theo}

The third equivalent statement was derived by Kolesnikov-Milman using the uniqueness of even-$L_p-$Minkowski problem.

\begin{theo}\label{t1.3} (\cite{KM}, page 66, Proposition 11.1)
 For each given dimension $n\geq2$ and $p\in[0,1)$, when restricted to origin-symmetric convex bodies in ${{\mathbb{R}}^{n+1}}$, the $L_p-$Brunn-Minkowski inequality \eqref{e1.3} is equivalent to the uniqueness of even-$L_p-$Minkowski problem
   \begin{equation}\label{e1.6}
     \det(\nabla^2u+uI)=fu^{p-1}, \ \ \forall x\in{\mathbb{S}}^n
   \end{equation}
on class of even functions.
\end{theo}

In order showing our main result, various local implications were also needed. As shown by Kolesnikov-Milman (\cite{KM}, page 67, Theorem 11.2), a local version of $p_*-BM-$inequality
     $$
       \frac{d^2}{d\varepsilon^2}\Bigg|_{\varepsilon=0}V(K+_p\varepsilon\cdot\varphi)^{\frac{p}{n}}\leq0, \ \ \forall\varphi\in C^{2,\alpha}({\mathbb{S}}^n)
     $$
implies local uniqueness of even-$L_p-$Minkowski problem for each $p>p_*$, where $K+_p\varepsilon\cdot\varphi$ is the convex body corresponding to the support function $(h_K^p+\varepsilon\varphi^p)^{\frac{1}{p}}$ for small $\varepsilon$. Combining with the implication of local $L_p-$Brunn-Minkowski inequality to local $p-BM-$inequality, one knows that the former implies local uniqueness of even-$L_p-$Minkowski problem in sense of Kolesnikov-Milman \cite{KM}.

Combining with Theorems \ref{t1.2} and \ref{t1.3}, our main Theorem \ref{t1.1} implies another three equivalent conjectures.

 \begin{theo}\label{t1.4}
  For each $p\in[0,1)$, both the $p-BM-$inequality and the uniqueness of even-$L_p-$Minkowski problem hold if and only if the third eigenvalue $\lambda_3(K)$ of Aleksandrov's problem \eqref{a.1} satisfies \eqref{a.2}. Moreover, the $L_p-$Minkowski inequality holds for $p\in(0,1)$ upon a same condition.
 \end{theo}

The contents of this paper are organized as follows: When localize on one of the convex bodies, we introduce a key necessary $p_*-$stable condition via second variation of the $L_p-$Brunn Minkowski inequality in Section 2. After a delicate comparison on localization concepts, we show that our $p_*-$stable condition implies the local uniqueness of even-$L_p-$Minkowski problem in Section 3. Utilizing the technics built by Chen-Huang-Li-Liu in \cite{CHLL}, the sufficiency of our $p_*$-stable condition implying the $L_p-$Brunn Minkowski inequality will be presented in Section 4. With the help of a new auxiliary variational scheme, we demonstrate the equivalence between the $p_*-$stable condition for $p_*=0$ and an inequality of the third eigenvalue of the Aleksandrov's problem in Section 5.

\vspace{30pt}

\section{Necessity of $p_*-$stable condition and second variation}

To initiate our arguments, we will construct a variation of the inequality by fixing one convex body and moving the other one. Considering first the case $p\in(0,1)$, we set
   $$
    h_\lambda\equiv\Big[(1-\lambda)h_K^p(x)+\lambda h_L^p(x)\Big]^{\frac{1}{p}}, \ \ h_L=h_K+\varepsilon\varphi
   $$
for even function $\varphi$. Noting that $h_\lambda$ is exactly the support function of Wulff body $W(h_\lambda)$ for small $\varepsilon$, the $L_p-$Brunn-Minkowski inequality is locally equivalent to state that the functional
  \begin{equation}\label{e2.1}
   I(\varepsilon)\equiv\Bigg(\int_{{\mathbb{S}}^n}h_\lambda\det(\nabla^2h_\lambda+h_\lambda I)\Bigg)\Bigg(\int_{{\mathbb{S}}^n}h_L\det(\nabla^2h_L+h_LI)\Bigg)^{-\lambda}
  \end{equation}
takes minimum at $\varepsilon=0$ for each given convex body $K$. Using the fact
  \begin{eqnarray}\nonumber\label{e2.2}
   &&\frac{d}{d\varepsilon}\Bigg|_{\varepsilon=0}h_\lambda=\lambda\varphi,\\
   &&\frac{d^2}{d\varepsilon^2}\Bigg|_{\varepsilon=0}h_\lambda=-\lambda(1-\lambda)(1-p)h_K^{-1}\varphi^2,
  \end{eqnarray}
there hold the vanishing first variation
  \begin{eqnarray*}
   &&\frac{d}{d\varepsilon}\Bigg|_{\varepsilon=0}I(\varepsilon)=(n+1)\lambda\Bigg(\int_{{\mathbb{S}}^n}\det(\nabla^2h_K+h_KI)\varphi\Bigg)\Bigg(\int_{{\mathbb{S}}^n}h_K\det(\nabla^2h_K+h_KI)\Bigg)^{-\lambda}\\
   &&-(n+1)\lambda\Bigg(\int_{{\mathbb{S}}^n}h_K\det(\nabla^2h_K+h_KI)\Bigg)^{-\lambda}\Bigg(\int_{{\mathbb{S}}^n}\det(\nabla^2h_K+h_KI)\varphi\Bigg)=0
  \end{eqnarray*}
and the second variation
  \begin{eqnarray*}
   &\displaystyle\frac{d^2}{d\varepsilon^2}\Bigg|_{\varepsilon=0}I(\varepsilon)=(n+1)^2\lambda(1-\lambda)\Bigg(\int_{{\mathbb{S}}^n}h_K\det(\nabla^2h_K+h_KI)\Bigg)^{-\lambda-1}\Bigg(\int_{{\mathbb{S}}^n}\det(\nabla^2h_K+h_KI)\varphi\Bigg)^2&\\
   &\displaystyle-(n+1)\lambda(1-\lambda)(1-p)\Bigg(\int_{{\mathbb{S}}^n}h_K\det(\nabla^2h_K+h_KI)\Bigg)^{-\lambda}\int_{{\mathbb{S}}^n}\det(\nabla^2h_K+h_KI)h_K^{-1}\varphi^2&\\
   &\displaystyle-(n+1)\lambda(1-\lambda)\Bigg(\int_{{\mathbb{S}}^n}h_K\det(\nabla^2h_K+h_KI)\Bigg)^{-\lambda}\Bigg(\int_{{\mathbb{S}}^n}U^{ij}(\varphi_{ij}+\varphi\delta_{ij})\varphi\Bigg),&
  \end{eqnarray*}
where $[U^{ij}]$ is the co-factor matrix of $[\nabla^2_{ij}h_K+h_K\delta_{ij}]$. When $p=0$, adapting the notations
   $$
    h_\lambda\equiv h^{\lambda}_L(x)h^{1-\lambda}_K(x), \ \ h_L=h_K+\varepsilon\varphi
   $$
and considering the functional $I(\cdot)$ defined in \eqref{e2.1}, one still have \eqref{e2.2}, the vanishing first variation and the formulated second variation above for $p=0$. So, we obtain the following necessary condition of $L_p-$Brunn-Minkowski inequality.

\begin{prop}\label{p2.1}
 Supposing that $L_p-$Brunn-Minkowski inequality \eqref{e1.3} was true for all $p\in(p_*,1), p_*\in[0,1)$, then ``$p_*-$stable condition"
   \begin{eqnarray}\nonumber\label{e2.3}
     &\displaystyle(n+1)\Bigg(\int_{{\mathbb{S}}^n} h_K\det(\nabla^2h_K+h_KI)\Bigg)^{-1}\Bigg(\int_{{\mathbb{S}}^n}\det(\nabla^2h_K+h_KI)\varphi\Bigg)^2&\\
     &\displaystyle\geq(1-p_*)\int_{{\mathbb{S}}^n}\det(\nabla^2h_K+h_KI)h_K^{-1}\varphi^2+\int_{{\mathbb{S}}^n}U^{ij}(\varphi_{ij}+\varphi\delta_{ij})\varphi&
   \end{eqnarray}
holds for all even functions $\varphi$ and $h_K$ satisfying
   \begin{equation}\label{e2.4}
    h_K>0, \ \ \nabla^2h_K+h_KI>0, \ \ \forall x\in{\mathbb{S}}^n.
   \end{equation}
\end{prop}

\vspace{30pt}

\section{Local uniqueness of even-$L_p-$Minkowski problem}

At the beginning of this section, let us start with a lemma of local $p-BM-$inequality in the following sense, which is different from that of Kolesnikov-Milman \cite{KM} by localization on $\lambda$.

\begin{lemm}\label{l3.1}
  Under the $p_*-$stable condition with some $p_*\in[0,1)$, given each $p\in(p_*,1)$ and $\varphi\in C^{2,\alpha}({\mathbb{S}}^n)$, the local $p-BM-$inequality holds in the sense of
     \begin{equation}\label{e3.1}
       V((1-\lambda)\cdot K+_p\lambda\cdot L)\geq\Big((1-\lambda)V^{\frac{p}{n}}(K)+\lambda V^{\frac{p}{n}}(L)\Big)^{\frac{n}{p}}, \ \ \forall \lambda\in[0,1],
     \end{equation}
 where $K,L$ are origin symmetric convex bodies satisfying
      \begin{equation}\label{e3.2}
          h_L=h_K+\varepsilon\varphi, \ \ \varepsilon\leq\varepsilon_0
      \end{equation}
 for some positive constant $\varepsilon_0$ depending only on $n, p, \varphi$ but not on $\lambda$.
\end{lemm}

\noindent\textbf{Proof.} When $p_*-$stable condition is true, by the calculations in Section 2, for any given $\varphi\in C^{2,\alpha}({\mathbb{S}}^n)$ and $\lambda\in[0,1]$,
   \begin{equation}\label{e3.3}
     \frac{d}{d\varepsilon}\Bigg|_{\varepsilon=0}I(\varepsilon)=0, \ \  \frac{d^2}{d\varepsilon^2}\Bigg|_{\varepsilon=0}I(\varepsilon)>0, \ \ h_L=h_K+\varepsilon\varphi
   \end{equation}
holds for any $p\in(p_*,1)$. We claim that there exists a positive constant $\varepsilon_0=\varepsilon_0(n, p, \varphi)$ such that
   \begin{equation}\label{e3.4}
     V((1-\lambda)\cdot K+_p\lambda\cdot L)\geq V^{1-\lambda}(K)V^\lambda(L), \ \ \forall\lambda\in[0,1]
   \end{equation}
for origin symmetric convex bodies $K,L$ satisfying \eqref{e3.2}. If not, for any $j\in{\mathbb{N}}$, there exist $\lambda_j\in[0,1]$ and $\varepsilon_j\to0$ such that
   \begin{equation}\label{e3.5}
     V((1-\lambda_j)\cdot K+_p\lambda_j\cdot L_j)<V^{1-\lambda_j}(K)V^{\lambda_j}(L_j)
   \end{equation}
holds for
   $$
    h_{L_j}=h_K+\varepsilon_j\varphi.
   $$
By Taylor's expansion, for each $j$, there exists $\xi_j\in(0,\varepsilon_j)$ such that
   \begin{equation}\label{e3.6}
     \frac{d^2}{d\varepsilon^2}\Bigg|_{\varepsilon=\xi_j}I_{\lambda_j}(\varepsilon)<0,
   \end{equation}
where
   $$
    I_{\lambda_j}(\varepsilon)\equiv\Bigg(\int_{{\mathbb{S}}^n}h_{\lambda_j}\det(\nabla^2h_{\lambda_j}+h_{\lambda_j}I)\Bigg)\Bigg(\int_{{\mathbb{S}}^n}h_{L_j}\det(\nabla^2h_{L_j}+h_{L_j}I)\Bigg)^{-\lambda_j}.
   $$
Assuming for a subsequence $\lambda_j\to\lambda_\infty\in[0,1]$, it yields from \eqref{e3.6} that,
   \begin{equation}\label{e3.7}
     \frac{d^2}{d\varepsilon^2}\Bigg|_{\varepsilon=0}I_{\lambda_\infty}(\varepsilon)\leq0
   \end{equation}
by sending $j\to\infty$, which contradicts with \eqref{e3.3}. The proof of claim was done. Now, replacing $K, L, \lambda$ in \eqref{e3.4} by
   $$
    \widetilde{K}\equiv\frac{K}{V^{\frac{1}{n}}(K)}, \ \ \widetilde{L}\equiv\frac{L}{V^{\frac{1}{n}}(L)}, \ \ \widetilde{\lambda}=\frac{\lambda V^{\frac{p}{n}}(L)}{(1-\lambda)V^{\frac{p}{n}}(K)+\lambda V^{\frac{p}{n}}(L)},
   $$
one derives the desired inequality \eqref{e3.1} upon \eqref{e3.2}. $\Box$\\

Next, we show that a local $p-BM$-inequality in sense of Kolesnikov-Milman \cite{KM} holds. Here some delicate comparisons are involved, since Kolesnikov-Milman \cite{KM} localize the inequalities in the parameter $\lambda$, but our arguments localize the inequalities in one of the convex body.

\begin{lemm}\label{l3.2}
  Supposing that $p_*-$stable condition holds for some $p_*\in[0,1)$, then for each $p\in(p_*,1)$ and $\varphi\in C^{2,\alpha}({\mathbb{S}}^n)$, the local $p-BM-$inequality holds in the sense of Kolesnikov-Milman
     \begin{equation}\label{e3.8}
       \frac{d^2}{d\varepsilon^2}\Bigg|_{\varepsilon=0}V(K+_p\varepsilon\cdot\varphi)^{\frac{p}{n}}\leq0,
     \end{equation}
  where $K+_p\varepsilon\cdot\varphi$ is the convex body corresponding to the support function $(h_K^p+\varepsilon\varphi^p)^{\frac{1}{p}}$ for small $\varepsilon$.
\end{lemm}

\noindent\textbf{Proof.} Setting
   $$
    F(h_\Omega)\equiv V^{\frac{p}{n}}(\Omega)
   $$
to be the functional of support function $h_\Omega$ of a convex body $\Omega$, it is inferred from \eqref{e3.1} (\cite{KM}, page 15, Lemma 3.1) that
   $$
    \frac{d^2}{d\lambda^2}V((1-\lambda)\cdot K+_p\lambda L)^{\frac{p}{n}}=\frac{d^2}{d\lambda^2}F\Big[((1-\lambda)h_K^p+\lambda h_L^p)^{\frac{1}{p}}\Big]\leq0
   $$
holds for all $\lambda\in[0,1]$ and given $K, L$ satisfying \eqref{e3.2}. Noting that
  \begin{eqnarray*}
   0&\geq&\frac{d^2}{d\lambda^2}F\Big[((1-\lambda)h_K^p+\lambda h_L^p)^{\frac{1}{p}}\Big]\\
   &=&\frac{1}{p^2}F''h_K^{2-2p}\Big(h_L^p-h_K^p\Big)^2+\frac{1}{p}\Bigg(\frac{1}{p}-1\Bigg)F'h_K^{1-2p}\Big(h_L^p-h_K^p\Big)^2.
  \end{eqnarray*}
After divided by $\varepsilon^2$, sending $\lambda\to0$ and then $\varepsilon\to0$, one obtains that
  \begin{equation}\label{e3.9}
    F''(h_K)\varphi^2+F'(h_K)(1-p)h_K^{-1}\varphi^2\leq0.
  \end{equation}
Another hand,
  \begin{eqnarray}\nonumber\label{e3.10}
   &&\frac{d^2}{d\varepsilon^2}\Bigg|_{\varepsilon=0}V(K+_p\varepsilon\cdot\varphi)^{\frac{p}{n}}=\frac{d^2}{d\varepsilon^2}\Bigg|_{\varepsilon=0}F\Big[(h_K^p+\varepsilon\varphi^p)^{\frac{1}{p}}\Big]\\
   &=&\frac{1}{p^2}F''(h_K)h_K^{4-2p}\varphi^2+\frac{1}{p}\Bigg(\frac{1}{p}-1\Bigg)F'(h_K)h_K^{3-4p}\varphi^2.
  \end{eqnarray}
So, \eqref{e3.8} follows from \eqref{e3.9} and \eqref{e3.10}. $\Box$\\

Combining with the implication result shown by Kolesnikov-Milman (\cite{KM}, page 67, Theorem 11.2), it yields from Lemma \ref{l3.2} the following local uniqueness of even-$L_p-$Minkowski problem.

\begin{prop}\label{p3.1}
  Supposing that $p_*-$stable condition holds for some $p_*\in[0,1)$, then for each $p\in(p_*,1)$ and original symmetric convex body $K$, there exists a $C^{2,\alpha}-$neighborhood $B_\sigma(K)$ of $K$ for each given $\alpha\in(0,1)$, such that the even-$L_p-$Minkowski problem admits a unique solution $K$ with respect to
    $$
     f\equiv h_K^{1-p}\det(\nabla^2h_K+h_KI), \ \ \forall x\in{\mathbb{S}}^n.
    $$
\end{prop}

\vspace{30pt}

\section{Sufficiency of $p_*-$stable condition}

In this section, we utilize the argument of Chen-Huang-Li-Liu in \cite{CHLL} to show that local uniqueness of even-$L_p-$Minkowski problem implies global uniqueness. Henceforth, we achieve the sufficiency of $p_*-$stable condition which implies the $L_p-$Brunn-Minkowski inequality. The simplified proofs are presented here for convenience of the readers. Different local to global argument can also be found a recent preprint by Milman \cite{M}. At first, let us quote a crucial {\it a-priori} lemma in \cite{CHLL}.

\begin{lemm}\label{l4.1}
  Suppose that $p\in[0,1]$ and $u$ is an even function satisfying
     \begin{equation}\label{e4.1}
       \det(\nabla^2u+uI)=fu^{p-1}, \ \ \forall x\in{\mathbb{S}}^n
     \end{equation}
  for some positive function $f\in C^\alpha({\mathbb{S}}^n)$. If $f$ satisfies that
     \begin{equation}\label{e4.2}
      C_1^{-1}\leq f\leq C_1,\ \ [f]_{C^\alpha({\mathbb{S}}^n)}\leq C_1, \ \ \forall x\in{\mathbb{S}}^n
     \end{equation}
  for some positive constant $C_1$, then
      \begin{equation}\label{e4.3}
        ||u||_{C^{2,\alpha}({\mathbb{S}}^n)}\leq C_2
      \end{equation}
  holds for some positive constant $C_2$ depending only on $n, p, C_1$ and $[f]_{C^\alpha({\mathbb{S}}^n)}$.
\end{lemm}

Since the eigenvalues of Beltrami-Laplace operator $\triangle=\triangle_{g_{{\mathbb{S}}^n}}$ on sphere are discrete, one can choose $p^*$ to be a number in $(0,1)$ such that $p^*-n-1$ is not an eigenvalue of $\triangle$. Now, for each $p\in(p_*,1)$ and even positive function $f\in C^\alpha({\mathbb{S}}^n)$, we set
   $$
     f_t\equiv tf+1-t, \ \ p_t\equiv tp+(1-t)p^*
   $$
and consider the functional equation
  \begin{equation}\label{e4.4}
    {\mathcal{F}}_t(u)\equiv\det(\nabla^2u+uI)-f_tu^{p_t-1}=0, \ \ \forall x\in{\mathbb{S}}^n
  \end{equation}
for $t\in[0,1]$. By Lemma \ref{l4.1}, there exists a positive constant $C_*$ independent of $t$, such that
   \begin{equation}\label{e4.5}
     C_*^{-1}\leq u\leq C_*, \ \ ||u||_{C^{2,\alpha}({\mathbb{S}}^n)}\leq C_*
   \end{equation}
holds for some $\alpha\in(0,1)$ and any zeros $u$ of ${\mathcal{F}}_t$. So, if one defines
  $$
   {\mathcal{O}}\equiv\Big\{u\in C^{2,\alpha}({\mathbb{S}}^n)\Big|\ (2C_*)^{-1}\leq u\leq 2C_*, \ \ ||u||_{C^{2,\alpha}({\mathbb{S}}^n)}\leq 2C_*\Big\},
  $$
it is clear that
   \begin{equation}\label{e4.6}
     {\mathcal{F}}_t^{-1}(0)\cap\partial{\mathcal{O}}=\emptyset.
   \end{equation}
Setting
  $$
   \Gamma\equiv\Big\{t\in[0,1]\Big|\ \mbox{there is a unique even positive classical solution to } \eqref{e4.4}\Big\},
  $$
 we will apply the method of topological degree to show the following result.

\begin{prop}\label{p4.1}
  Supposing that $p_*-$stable condition holds, then for each $p\in(p_*,1)$ and even positive function $f\in C^\alpha({\mathbb{S}}^n), \alpha\in(0,1)$,
    \begin{equation}
      \Gamma=[0,1].
    \end{equation}
\end{prop}

\noindent\textbf{Proof.} At first, it is inferred from a uniqueness result by Brendle-Choi-Daskalopoulos \cite{BCD} that $u_0\equiv1$ is the unique solution to \eqref{e4.4} when $t=0$. Moreover, the linearized equation of \eqref{e4.4} at $u_0$ is given by
  \begin{equation}\label{e4.8}
    \triangle\varphi+(n+1-p^*)\varphi=0, \ \ \mbox{ on } {\mathbb{S}}^n.
  \end{equation}
Using the choice of $p^*$, there is only trivial zero solution of \eqref{e4.8}. So, we have
   $$
    deg({\mathcal{F}}_0,{\mathcal{O}},0)=1\not=0.
   $$
Combining with the preservation property of topological degree \cite{L,N}, there holds
   \begin{equation}\label{e4.9}
     deg({\mathcal{F}}_t,{\mathcal{O}},0)=1\not=0, \ \ \forall t\in[0,1].
   \end{equation}
So, the solvability of \eqref{e4.4} was done for all $t\in[0,1]$. It remains to show that the solution of \eqref{e4.4} is unique for each $t\in[0,1]$. We follow the argument in \cite{CHLL} to show that $\Gamma$ is both closed and semi-open in the sense of
    $$
     t_0\in\Gamma\Rightarrow [t_0,t_0+\varepsilon)\in\Gamma \mbox{ for some } \varepsilon>0.
    $$
In fact, supposing that $t_0\in\Gamma$, we claim that for some $\varepsilon>0$, the solution of \eqref{e4.4} for $t\in(t_0,t_0+\varepsilon)$ is unique. If not, for each $j\in{\mathbb{N}}$, there exists at least two different solutions $u_j$ and $v_j$ to \eqref{e4.4} with $t=t_0+\varepsilon_j$, where $\varepsilon_j\to0$. Let $u_0$ a solution of \eqref{e4.4} at $t=t_0$. By {\it a-priori} Lemma \ref{l4.1}, there must be
    \begin{equation}\label{e4.10}
      \lim_{j\to\infty}u_j=u_0, \ \ \lim_{j\to\infty}v_j=u_0, \ \mbox{ uniformly on } C^{2,\alpha'}({\mathbb{S}}^n)
    \end{equation}
for some $\alpha'\in(0,\alpha)$. Otherwise, it will yield two solutions of \eqref{e4.4} at $t=t_0$, which contradicts with local uniqueness result at $t=t_0$. Even though, \eqref{e4.10} still contradicts with our Proposition \ref{p3.1} at $t=t_0+\varepsilon_j$ for large $j$. So, our claim of semi-openness of $\Gamma$ holds true.

Now, let us show $\Gamma$ is also closed by following the argument in \cite{CHLL}. Supposing on the contrary, there exists $t_0\in(0,1)$ such that the solution of $\eqref{e4.4}$ is unique for $t\in(t_0-\varepsilon_0, t_0), \varepsilon_0>0$ but not for $t=t_0$. Let $h_1, h_2$ be two different solutions of \eqref{e4.4} at $t=t_0$. By lemma \ref{l4.2} under below, when $t\in(0,t_0)$ closing to $t_0$, there exist two solutions $h_1', h'_2$ of \eqref{e4.4} which are close to $h_1, h_2$ respectively. This contradicts with our uniqueness assumption of \eqref{e4.4} for $t\in(t_0-\varepsilon_0, t_0)$. So, desired conclusion of Proposition \ref{p4.1} follows. $\Box$\\

Next, let us complete the proof of Proposition \ref{p4.1} by proving the following lemma.

\begin{lemm}\label{l4.2}
  Suppose \eqref{e4.4} admits an even positive classical solution $h_{t_0}$ for some $t_0\in(0,1)$ and even positive $f_{t_0}\in C^\alpha({\mathbb{S}}^n)$. For any $\delta>0$, there exists $\varepsilon>0$ such that \eqref{e4.4} admits an even positive classical solution $h_t$ satisfying
     \begin{equation}\label{e4.11}
       ||h_t-h_{t_0}||\leq\delta, \ \ t_0-\varepsilon<t<t_0.
     \end{equation}
\end{lemm}

\noindent\textbf{Proof.} At first, noting that the solution $h_{t_0}$ is also a solution of
   \begin{equation}\label{e4.12}
     \det(\nabla^2h+hI)=\overline{f}_0h^{p_t-1}, \ \ \overline{f}_0\equiv f_{t_0}h_{t_0}^{p_{t_0}-p_t},
   \end{equation}
whose linearized equation is given by
   $$
   U^{ij}(\phi_{ij}+\phi\delta_{ij})=(p_t-1)\frac{\phi}{h_{t_0}}\det(\nabla^2{h_{t_0}}+{h_{t_0}}I)
   $$
for $[U^{ij}]$ being the co-factor matrix of $[\nabla^2h_{t_0}+h_{t_0}I]$. Because the eigenvalues of linear operator
   $$
    {\mathcal{L}}\phi\equiv\frac{h_{t_0}}{\det(\nabla^2{h_{t_0}}+{h_{t_0}}I)}U^{ij}(\phi_{ij}+\phi\delta_{ij})
   $$
are discrete, there exists $\varepsilon_0$ small enough such that the kernel of ${\mathcal{L}}-(p_t-1)I$ is a null space for any $t_0-\varepsilon_0<t<t_0$. Consider the functional
  $$
   {\mathcal{G}}_\lambda(u)\equiv\det(\nabla^2u+uI)-\overline{f}_\lambda u^{p_t-1}, \ \ \overline{f}_\lambda\equiv(1-\lambda)\overline{f_0}+\lambda f_t.
  $$
Using the fact
   \begin{equation}\label{e4.13}
    ||\overline{f}_\lambda-\overline{f}_0||_{C^{\alpha}({\mathbb{S}}^n)}=\lambda||f_t-f_{t_0}h_{t_0}^{p_{t_0}-p_t}||_{C^{\alpha}({\mathbb{S}}^n)}\leq\varepsilon_1\ll1
   \end{equation}
for any $\lambda\in[0,1]$ and small $\varepsilon\leq\varepsilon_0$, by local uniqueness Proposition \ref{p3.1} and {\it a-priori} Lemma \ref{l4.1}, there exists a small positive constant $\varepsilon_2<\sigma/4$ such that
  \begin{equation}\label{e4.14}
   \mbox{either } ||h_\lambda-h_0||_{C^{2,\alpha}({\mathbb{S}}^n)}>\sigma, \ \mbox{ or } ||h_\lambda-h_0||_{C^{2,\alpha}({\mathbb{S}}^n)}<\varepsilon_2,
  \end{equation}
where $h_\lambda$ is a possible zero of ${\mathcal{G}}_\lambda(\cdot)$ and $\sigma$ is a positive constant coming from Proposition \ref{p3.1}. So, if one considers zeros of ${\mathcal{G}}_\lambda(\cdot)$ on bounded open set
   $$
    {\mathcal{W}}\equiv\Big\{u\in C^{2,\alpha}({\mathbb{S}}^n)\Big|\ ||u-h_0||_{C^{2,\alpha}({\mathbb{S}}^n)}<2\varepsilon_2\Big\},
   $$
there holds
  \begin{equation}\label{e4.15}
    {\mathcal{G}}^{-1}_\lambda(0)\cap\partial{\mathcal{W}}=\emptyset.
  \end{equation}
Moreover, $h_0=h_{t_0}$ is a zero of ${\mathcal{G}}_0(\cdot)$, whose linearized equation has null kernel when $\varepsilon\leq\varepsilon_0$. We thus conclude that
  \begin{equation}\label{e4.16}
   deg({\mathcal{G}}_\lambda,{\mathcal{W}},0)=deg({\mathcal{G}}_0,{\mathcal{W}},0)=1\not=0
  \end{equation}
for each $\lambda\in[0,1]$. So, ${\mathcal{G}}_1(\cdot)$ has a desired zero inside ${\mathcal{W}}$. The proof of Lemma \ref{l4.2} was done. $\Box$\\

Let us sum the following result at the end of this section.

\begin{theo}\label{t4.1}
   Supposing that $p_*-$stable condition is true for some $p_*\in[0,1)$, then for each $p\in(p_*,1)$, the global $L_p-$Brunn-Minkowski inequality holds.
\end{theo}

\noindent\textbf{Proof.} The conclusion of this theorem is exactly a corollary of Proposition \ref{p4.1} and Theorem \ref{t1.3}. $\Box$\\

\vspace{30pt}

\section{$p_*-$stable condition \eqref{e2.3} and third eigenvalue of Aleksandrov problem}

Now, we shall reduce the $0-$stable condition to a variational scheme as follows. Consider
 \begin{eqnarray*}
   &\displaystyle\inf_{\varphi\in{\mathcal{A}}}J(\varphi), \ \ {\mathcal{A}}\equiv\Bigg\{\varphi\in H^1({\mathbb{S}}^n)\Bigg|\ \int_{{\mathbb{S}}^n}h_K^{-1}\det(\nabla^2h_K+h_KI)\varphi^2=1, \ \varphi\mbox{ even}\Bigg\}&\\
   &\displaystyle J(\varphi)\equiv-\int_{{\mathbb{S}}^n}U^{ij}(\nabla^2_{ij}\varphi+\varphi\delta_{ij})\varphi&\\
   &\displaystyle+(n+1)\Bigg(\int_{{\mathbb{S}}^n}h_K\det(\nabla^2h_K+h_KI)\Bigg)^{-1}\Bigg(\int_{{\mathbb{S}}^n}\det(\nabla^2h_K+h_KI)\varphi\Bigg)^2&
 \end{eqnarray*}
for fixed convex body $K$ and $[U^{ij}]=co[\nabla^2_{ij}h_K+h_K\delta_{ij}]$. The Euler-Lagrange equation of the variational problem is given by
   \begin{equation}\label{e5.1}
    -U^{ij}(\nabla^2_{ij}\varphi+\varphi\delta_{ij})+\lambda_1(\varphi)\det(\nabla^2h_K+h_KI)=\lambda_2h_K^{-1}\det(\nabla^2h_K+h_KI)\varphi
   \end{equation}
for some positive constants $\lambda_2$ and
   $$
    \lambda_1(\varphi)\equiv(n+1)\Bigg(\int_{{\mathbb{S}}^n}h_K\det(\nabla^2h_K+h_KI)\Bigg)^{-1}\Bigg(\int_{{\mathbb{S}}^n}\det(\nabla^2h_K+h_KI)\varphi\Bigg)
   $$

We have the following equivalence.

\begin{prop}\label{p5.1}
 Fixing an origin-symmetric convex body $K\subset{\mathbb{R}}^{n+1}$, the $p_*-$stable condition \eqref{e2.3} holds with $p_*=0$ if and only if the identity
     \begin{equation}\label{e5.2}
      \inf_{\varphi\in{\mathcal{A}}}J(\varphi)=1
     \end{equation}
  holds for the same convex body $K$.
\end{prop}

Let us recall a result of Aleksandrov type (\cite{A2}, section 6) concerning the first and second eigenvalues of the eigenvalue problem
   \begin{equation}\label{e5.3}
      -U^{ij}(\nabla^2_{ij}\phi+\phi\delta_{ij})=\lambda h_K^{-1}\det(\nabla^2h_K+h_KI)\phi.
   \end{equation}

\begin{prop}\label{p5.2}
  Considering the eigenvalue problem \eqref{e5.3}, the eigenfunctions corresponding to the zero eigenvalue $\lambda=0$ are given by $\phi=a\cdot x$ for some vector $a\in{\mathbb{R}}^{n+1}$. Moreover, $\lambda=-n$ is the unique negative eigenvalue of \eqref{e5.3}, whose eigenfunctions can only by $\kappa h_K$ for each $\kappa\not=0$.
\end{prop}

We have the following equivalence between the $p_*-$stable condition and the third eigenvalue of the Aleksandrov's eigenvalue problem \eqref{e5.3}.

\begin{theo}\label{t5.1}
 Fixing an origin-symmetric convex body $K\subset{\mathbb{R}}^{n+1}$, the $p_*-$stable condition \eqref{e2.3} holds for $p_*=0$ if and only if the third eigenvalue of \eqref{e5.3} satisfies that
     \begin{equation}\label{e5.4}
       \lambda_3(K)\geq1
     \end{equation}
  for the same convex body $K$.
\end{theo}

We need the following inequality.

\begin{lemm}\label{l5.1}
  Given origin-symmetric convex body $K$, the inequality
    \begin{eqnarray}\label{e5.5}
      &\int_{{\mathbb{S}}^n}\varphi U^{ij}(\nabla^2_{ij}\varphi+\varphi\delta_{ij})+\lambda_3(K)\int_{{\mathbb{S}}^n}h_K^{-1}\det(\nabla^2h_K+h_KI)\varphi^2&\\ \nonumber
      &\leq(n+\lambda_3(K))\Big(\int_{{\mathbb{S}}^n}h_K\det(\nabla^2h_K+h_KI)\Big)^{-1}\Big(\int_{{\mathbb{S}}^n}\det(\nabla^2h_K+h_KI)\varphi\Big)^2&
    \end{eqnarray}
  holds for any even function $\varphi$ on ${\mathbb{S}}^n$, where $\lambda_3(K)$ is the third eigenvalue of \eqref{e5.3}. Moreover, equality holds if and only if
     \begin{equation}\label{e5.6}
      \phi\equiv\varphi-\frac{\int_{{\mathbb{S}}^n}\det(\nabla^2h_K+h_KI)\varphi}{\int_{{\mathbb{S}}^n}h_K\det(\nabla^2h_K+h_KI)}h_K
     \end{equation}
  is an eigenfunction corresponding to the third eigenvalue of \eqref{e5.3} or identical to zero.
\end{lemm}

\noindent\textbf{Proof.} Letting $\phi$ to be defined in \eqref{e5.6}
and noting that $a\cdot x$ are odd functions for all vectors $a\in{\mathbb{R}}^{n+1}$, we have $\phi$ is orthonormal to the first and second eigen-spaces of \eqref{e5.3}. Hence, using the fact
    $$
     -\int_{{\mathbb{S}}^n}\phi U^{ij}(\nabla^2_{ij}\phi+\phi\delta_{ij})\geq\lambda_3(K)\int_{{\mathbb{S}}^n}h_K^{-1}\det(\nabla^2h_K+h_KI)\phi^2
    $$
and writing back to $\varphi$, we derived the desired inequality \eqref{e5.5}. $\Box$\\

Now, one can complete the proof of Theorem \ref{t5.1}.

\noindent\textbf{Proof of Theorem \ref{t5.1}.} Let us first show the sufficiency of \eqref{e5.4}. For each even function $\varphi\in{\mathcal{A}}$, we have
   \begin{equation}\label{e5.7}
     \int_{{\mathbb{S}}^n}h_K^{-1}\det(\nabla^2h_K+h_KI)\varphi^2=1.
   \end{equation}
Thus, it follows from \eqref{e5.5} and \eqref{e5.7} that
   \begin{eqnarray}\label{e5.8}
     &J(\varphi)\geq1+(\lambda_3(K)-1)\int_{{\mathbb{S}}^n}h_K^{-1}\det(\nabla^2h_K+h_KI)\varphi^2&\\ \nonumber
     &-(\lambda_3(K)-1)\Bigg(\int_{{\mathbb{S}}^n}h_K\det(\nabla^2h_K+h_KI)\Bigg)^{-1}\Bigg(\int_{{\mathbb{S}}^n}\det(\nabla^2h_K+h_KI)\varphi\Bigg)^2.&
   \end{eqnarray}
Using the H\"{o}lder inequality
   $$
    \int_{{\mathbb{S}}^n}\det(\nabla^2h_K+h_KI)\varphi\leq\sqrt{\int_{{\mathbb{S}}^n}h_K^{-1}\det(\nabla^2h_K+h_KI)\varphi^2}\sqrt{\int_{{\mathbb{S}}^n}h_K\det(\nabla^2h_K+h_KI)},
   $$
if \eqref{e5.4} holds, one has $J(\varphi)\geq1$. Combining with the fact
   $$
    J\Bigg(\frac{h_K}{\sqrt{\int_{{\mathbb{S}}^n}h_K\det(\nabla^2h_K+h_KI)}}\Bigg)=1,
   $$
the sufficiency of \eqref{e5.4} follows by Proposition \ref{p5.1}. It remains to show the necessity of \eqref{e5.4}. Actually, if not, there exists at least one eigenfunction $\phi$ corresponding to the third eigenvalue $\lambda_3(K)<1$. Setting
  $$
    \varphi_\varepsilon\equiv \frac{h_K+\varepsilon\phi}{\sqrt{\int_{{\mathbb{S}}^n}h_K^{-1}\det(\nabla^2h_K+h_KI)({h_K+\varepsilon\phi})^2}}
  $$
for a given small constant $\varepsilon>0$, the above calculation shows that $J(\varphi_\varepsilon)<1$, which contradicts with the validity of $p_*-$stable condition for $p_*=0$. The proof was done. $\Box$\\

\vspace{30pt}

\section*{Acknowledgments}

The author would like to express his deepest gratitude to Professors Xi-Ping Zhu, Kai-Seng Chou, Xu-Jia Wang and Neil Trudinger for their constant encouragements and warm-hearted helps. Special thanks were owed to Professors Shi-Bing Chen and Emanuel Milman for some valuable conversations and comments. This paper was also dedicated to the memory of Professor Dong-Gao Deng.\\

\end{document}